\documentclass{amsart}
\usepackage{graphicx}
\usepackage{amssymb,amscd,amsthm,amsxtra}
\usepackage{latexsym}

\newtheorem{thm}{Theorem}[section]
\newtheorem{cor}[thm]{Corollary}
\newtheorem{lem}[thm]{Lemma}
\newtheorem{prop}[thm]{Proposition}
\theoremstyle{definition}
\newtheorem{defn}[thm]{Definition}
\theoremstyle{remark}
\newtheorem{rem}[thm]{Remark}

\numberwithin{equation}{section}

\begin{document}

\title[Singular minimizers]{Some singular minimizers in low dimensions in the calculus of variations}
\author{Connor Mooney}
\author{Ovidiu Savin}
\address{Department of Mathematics, Columbia University, New York, NY 10027}
\email{\tt  cmooney@math.columbia.edu}
\email{\tt savin@math.columbia.edu}

\begin{abstract}
 We construct a singular minimizing map $\bf u $ from $\mathbb{R}^3$ to $\mathbb{R}^2$ of a smooth uniformly convex functional of the form $\int_{B_1} F(D{\bf u})\,dx$.
\end{abstract}
\maketitle
\section{Introduction}
In this paper we consider minimizers of functionals of the form
\begin{equation}\label{Functional}
 \int_{B_1} F(D{\bf u})\,dx
\end{equation}
where ${\bf u} \in H^1(B_1)$ is a map from $\mathbb{R}^n$ to $\mathbb{R}^m$ and $F$ is a smooth, uniformly convex function on $M^{m \times n}$ with bounded
second derivatives. By a minimizer we understand a map $\bf u$ for which the integral above increases after we perform any smooth deformation of $\bf u$, with compact support in $B_1$. If $F$ satisfies these conditions then minimizers are unique subject to their own boundary condition. Moreover $\bf u$ is a minimizer if and only if it solves the Euler-Lagrange system 
\begin{equation}\label{MinimizerEquation}
 \text{div}(\nabla F(D{\bf u})) = 0,
\end{equation} 
in the sense of distributions.

The regularity of minimizers of \eqref{Functional} is a well-studied problem. Morrey \cite {Mo} showed that in dimension $n = 2$
all minimizers are smooth. This is also true in the scalar case $m = 1$ by the classical results of De Giorgi and Nash \cite{DG1},\cite{Na}. In the scalar case, the regularity is obtained by differentiating equation \eqref{MinimizerEquation} and treating the problem as a linear equation with bounded measurable coefficients.
An example of De Giorgi \cite{DG2} shows that these techniques cannot be extended to the case $m \geq 2$. Another example due to Giusti and Miranda \cite{GM2} shows that elliptic systems do not have regularity even when the coefficients depend only on ${\bf u}$.
On the other hand it is known that minimizers of \eqref{Functional} are smooth away from a closed singular set of Hausdorff $n-p$ dimensional measure zero for
some $p > 2$, see \cite{GM1}, \cite{GG}. (In fact, if $F$ is uniformly quasi-convex then minimizers are smooth away from a closed set of Lebesgue
measure zero, see Evans \cite{E2}). However, the singular set may be non-empty. We will discuss some interesting examples below.

The main result of this paper is a counterexample to the regularity of minimizers of \eqref{Functional} when $n = 3$ and $m = 2$, which are the
optimal dimensions in light of the previous results. The existence of such minimizing maps from $\mathbb{R}^3$ to $\mathbb{R}^3$ or from $\mathbb{R}^3$ to $\mathbb{R}^2$ is stated as an open problem in the book of Giaquinta (see \cite{Gi}, p. 61).

The first example of a singular minimizer of \eqref{Functional} is due to Ne\v{c}as \cite{Ne}. He considered the homogeneous degree one map 
$${\bf u}(x)= \frac{x \otimes x}{|x|}$$ 
from $\mathbb{R}^n$ to $\mathbb{R}^{n^2}$
for $n$ large, and constructed explicitly a smooth uniformly convex $F$ on $M^{n^2 \times n}$ for which ${\bf u}$ minimizes \eqref{Functional}.
Later Hao, Leonardi and Ne\v{c}as \cite{HLN} improved the dimension to $n = 5$ using
\begin{equation}\label{SYExample}
{\bf u}(x) = \frac{x \otimes x}{|x|} - \frac{|x|}{n}I.
\end{equation}

The values of \eqref{SYExample} are symmetric and traceless, and thus lie in a $n(n+1)/2 - 1$ dimensional subspace of $M^{n \times n}$.
\v{S}ver\'{a}k and Yan \cite{SY1} showed that the map \eqref{SYExample} is a counterexample for $n=3,\,m=5$. Their approach
is to construct a quadratic null Lagrangian $L$ which respects the symmetries of ${\bf u}$, such that $\nabla L = \nabla F$ on $D{\bf u}(B_1)$ for 
some smooth, uniformly convex $F$ on $M^{5 \times 3}$. The Euler-Lagrange system $\text{div}(\nabla F(D{\bf u})) = \text{div}(\nabla L(D{\bf u})) = 0$
then holds automatically. In \cite{SY2} they use the same technique to construct a non-Lipschitz minimizer with $n = 4,\,m=3$ coming from the Hopf fibration. 
To our knowledge, these are the lowest-dimensional examples to date.

Our strategy is different and it is based on constructing a homogenous of degree one minimizer in the scalar case for an integrand which is convex but has ``flat pieces".

An interesting problem about the regularity of minimizers occurs in the scalar case when considering in \eqref{Functional} convex integrands $F:\mathbb{R}^n \to \mathbb{R}$ for which the uniform convexity of $F$ fails on some compact set $\mathcal S$. Assume for simplicity that $F$ is smooth outside the degeneracy set $\mathcal S$, and also that $F$ satisfies the usual quadratic growth at infinity. One key question is whether or not the gradient $\nabla u$ localizes as we focus closer and closer to a point 
$x_0 \in B_1$. In \cite{DS} it was proved that, in dimension $n=2$, the sets $\nabla u(B_\varepsilon(x_0))$
decrease uniformly as $\varepsilon \to 0$ either to a point outside $\mathcal S$, or to a connected subset of $\mathcal S$. 
In Theorem \ref{ScalarExample} below we show that this ``continuity property" of $\nabla u$ does not hold in dimension $n=3$ when the set $\mathcal S$ is the union of two disconnected convex sets. 
 We remark that, as in the $p$-Laplace equation, it is relatively standard (see \cite{E1,CF}) to obtain the continuity of $\nabla u$ outside the convex hull $\mathcal S^c$ of $\mathcal S$.

Let $w$ be the homogeneous degree one function
$$w(x_1,x_2) = \frac{x_2^2 - x_1^2}{\sqrt{2(x_1^2 + x_2^2)}}=\frac{-1}{\sqrt 2} \, \, r  \, \cos 2 \theta,$$
and let $u_0$ be the function on $\mathbb{R}^3$ obtained by revolving $w$ around the $x_1$ axis,
$$u_0(x_1,x_2,x_3) = w\left(x_1,\sqrt{x_2^2 + x_3^2}\right).$$
We show that $u_0$ solves a degenerate elliptic equation that is uniformly elliptic away from the cone 
$$K_0 = \{x_1^2 > x_2^2+x_3^2\}.$$

\begin{thm}\label{ScalarExample}
 For any $\delta > 0$ there exists a convex function $G_0 \in C^{1,1-\delta} (\mathbb{R}^3)$ which is linear on two bounded convex sets containing $\nabla u_0(K_0)$,
 uniformly convex and smooth away from these two convex sets, such that $u_0$ is a minimizer of the functional $$\int_{B_1} G_0(\nabla u_0)\,dx.$$
\end{thm}

We use $u_0$ and $G_0$ to construct a singular minimizing map from $\mathbb{R}^3$ to $\mathbb{R}^2$. Rescaling $u_0$ we obtain a function $u^1$ that solves an equation
that is uniformly elliptic away from a thin cone around the $x_1$ axis, and switching the $x_1$ and $x_3$ axes we get an analogous function
$u^2$. Then ${\bf u} = (u^1,u^2)$ is a minimizing map for
 $$F_0(p^1,p^2) := G_1(p^1) + G_2(p^2),$$ which is a convex function defined on $\mathbb{R}^6 \cong M^{2 \times 3}$. Notice that the Euler-Lagrange system
$\text{div}(\nabla F_0(D{\bf u}))=0$ is de-coupled, and $F_0$ fails to be uniformly convex or smooth in certain regions. However, a key observation is that $F_0$ separates quadratically from its tangent planes when restricted to the image of $D {\bf u}$. We obtain our example by making a small perturbation of $F_0$.

More specifically, let
$$u^1(x_1,x_2,x_3) = u_0(x_1/2,x_2,x_3), \quad u^2(x_1,x_2,x_3) = u^1(x_3,x_2,x_1)$$
and let 
\begin{equation}\label{MainExample}
 {\bf u} = (u^1,u^2).
\end{equation}

Our main theorem is:

\begin{thm}\label{main}
The map \eqref{MainExample} is a minimizer of $$\int_{B_1} F(D{\bf u})\,dx$$ for some smooth, uniformly convex $F : M^{2 \times 3} \rightarrow \mathbb{R}$.
\end{thm}

The paper is organized as follows. In Section \ref{MainProof} we state a convex extension lemma and the key proposition, which asserts the existence
of a suitable smooth small perturbation of $G_0$. We then use them to prove Theorem \ref{main}.
In Section \ref{Constructions} we prove the key proposition. This section contains most of the technical details.
In Section \ref{Appendix} we prove the extension lemma and some technical inequalities needed for the key proposition. Finally,
at the end of Section \ref{Appendix} we outline how to prove Theorem \ref{ScalarExample}.

\section{Key Proposition and Proof of Theorem \ref{main}}\label{MainProof}

In this section we state the extension lemma and the key proposition. We then use them to prove Theorem \ref{main}.

The function $F_0$ defined in the Introduction is not uniformly convex in $M^{2 \times 3}$, but it
separates quadratically from its tangent planes on the image of $D \bf u$ which, by the one-homogeneity of $\bf u$, is the two dimensional surface $D{\bf u}(S^2)$. The quadratic separation holds on this surface since $G_1$ is uniformly convex in the region where $G_2$ is flat and vice versa.
We would like to find a uniformly convex extension of $F_0$ with the same tangent planes on $D{\bf u}( \partial B_1)$.

\subsection{Extension Lemma}

The extension lemma gives a simple criterion for deciding when the tangent planes on a smooth surface can be extended to a global smooth, uniformly convex function. 
Let $\Sigma$ be a smooth compact, embedded surface in $\mathbb{R}^n$ of any dimension.

\begin{lem}\label{ExtensionLemma}
Let $G$ be a smooth function and ${\bf v}$ a smooth vector field on $\Sigma$ such that
\begin{equation}\label{QuadSepCondition}
 G(y) - G(x) - {\bf v}(x) \cdot (y-x) \ge \gamma |y-x|^2,
\end{equation}
for any $x, \,y \in \Sigma$ and some $\gamma > 0$. Then there exists a global smooth function $F$ such that $F = G$ and $\nabla F = {\bf v}$ on $\Sigma$, 
and $D^2F \ge \gamma I$.
\end{lem}

The idea of the proof is to first make a local extension by adding a large multiple of the square of distance from $\Sigma$. We then make an extension to all of $\mathbb{R}^n$ by taking the supremum of tangent paraboloids to the local extension.
Finally we mollify and glue the local and global extensions. We postpone the proof to the appendix, Section \ref{Appendix}. We also
record an obvious corollary.

\begin{defn}\label{SeparationDefinition}
 Let $G$ be a smooth function on an open subset $O$ of $\mathbb{R}^n$. We define the separation function $S_G$ on $O \times O$ by
 $$S_{G}(x,y) = G(y) - G(x) - \nabla G(x) \cdot (y-x).$$
\end{defn}

\begin{cor}\label{ExtensionCorollary}
Assume that $G$ is a smooth function in a neighborhood of $\Sigma$ such that $S_G(x,y) \ge \gamma |y-x|^2$
for any $x, \,y \in \Sigma$ and some $\gamma > 0$. Then there exists a global smooth, uniformly convex function $F$ such that 
$F = G$ and $\nabla F = \nabla G$ on $\Sigma$.
\end{cor}
\subsection{Key Proposition}
In this section we state the key proposition. We first give the setup for the statement. Recall that $w = (x_2^2-x_1^2)/\sqrt{2(x_1^2+x_2^2)}$. Let
$$\Gamma = \nabla w(B_1 - \{0\}) = \nabla w(S^1).$$
We describe $\Gamma$ as a collection of four congruent curves. The part of $\Gamma$ in the region $\{p_2 \geq |p_1|\}$ can be written as a graph
$$\Gamma_1 = \{(p_1, \varphi(p_1))\}$$
for $p_1 \in [-1,1]$, where $\varphi$ is even, uniformly convex, tangent to $p_1^2=p_2^2$ at $\pm 1$, and
separates from these lines like $(\text{dist})^{3/2}$. We will give a more precise description of $\varphi$ in Section \ref{Constructions}. 

\begin{figure}
 \centering
    \includegraphics[scale=0.35]{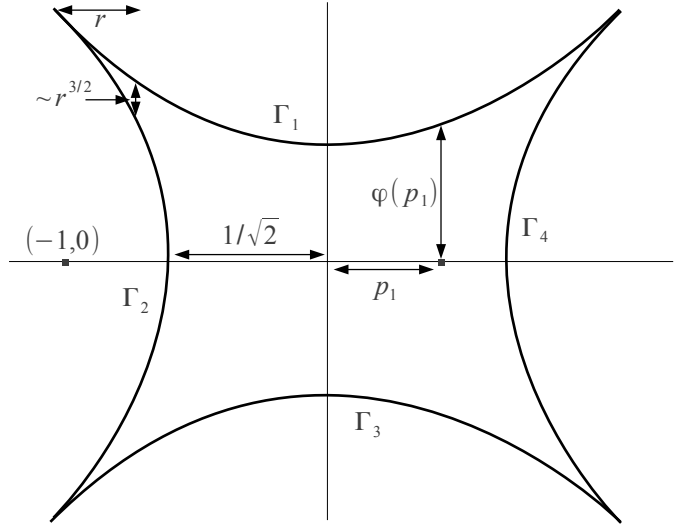}
 \caption{$\Gamma$ consists of four identical curves separating from the lines $p_1^2=p_2^2$ like $\text{dist}^{3/2}.$}
 \label{Gammapic}
\end{figure}

The other pieces of $\Gamma$ can be written
$$\Gamma_2 = \{-\varphi(p_2),p_2\}, \quad \Gamma_3 = \{p_1, -\varphi(p_1)\}, \quad \Gamma_4 = \{\varphi(p_2),p_2\}$$
for $p_i \in [-1,1]$, representing the left, bottom and right  pieces of $\Gamma$ (see figure \ref{Gammapic}).

Recall that $u_0 = w\left(x_1,\sqrt{x_2^2+x_3^2}\right)$. Then
$$\Omega = \nabla u_0(S^2)$$
is the surface obtained by revolving $\Gamma$ around the $p_1$ axis. Let $\Omega_R \subset \Omega$ be the surface obtained by revolving
$\Gamma_1$ around the $p_1$ axis. 

In the statement below, $\delta$ and $\gamma$ are small positive constants depending on $\varphi$.

\begin{prop}\label{KeyProposition}
For any $\epsilon > 0$ there exists a smooth function $G$ defined in a neighborhood of $\Omega$ such that 
$$\text{div}(\nabla G(\nabla u_0)) = 0 \quad \quad \mbox{in} \quad B_1 \setminus\{0\},$$ and
\begin{enumerate}
\item If $p \in \Omega_R \cap \{-1 + \delta \leq p_1 \leq 1-\delta\}$ then $S_{G}(p,q) \ge \gamma|p-q|^2$
for all $q \in \Omega$, 
\item $S_{G}(p,q) \geq -\epsilon |p-q|^2$ otherwise for $p,\,q \in \Omega$.
\end{enumerate}
\end{prop}

We delay the proof of this proposition to Section \ref{Constructions}, and use it now to prove Theorem \ref{main}.

\subsection{Proof of Theorem \ref{main}}
Recall that
$$u^1(x_1,x_2,x_3) = u_0(x_1/2,x_2,x_3), \quad u^2(x_1,x_2,x_3) = u^1(x_3,x_2,x_1),$$
and let 
$$G_1(p_1,p_2,p_3) = G(2p_1,p_2,p_3), \quad G_2(p_1,p_2,p_3) = G_1(p_3,p_2,p_1).$$
Then by Proposition \ref{KeyProposition} we have $\text{div}(\nabla G_i(\nabla u^i)) = 0$. Let 
$$\Sigma = D{\bf u}(B_1).$$
Since $D^2u^1$ has rank $2$ away from the cone 
$$K_1 = \left\{x_1^2 \geq 4(x_2^2+x_3^2)\right\}$$ and similarly $D^2u^2$ has rank $2$ away from
$$K_2 = \left\{x_3^2 \geq 4(x_1^2 + x_2^2)\right\},$$ 
it is easy to see that $\Sigma$ is a smooth embedded surface in $\mathbb{R}^6$.

\begin{figure}
 \centering
    \includegraphics[scale=0.45]{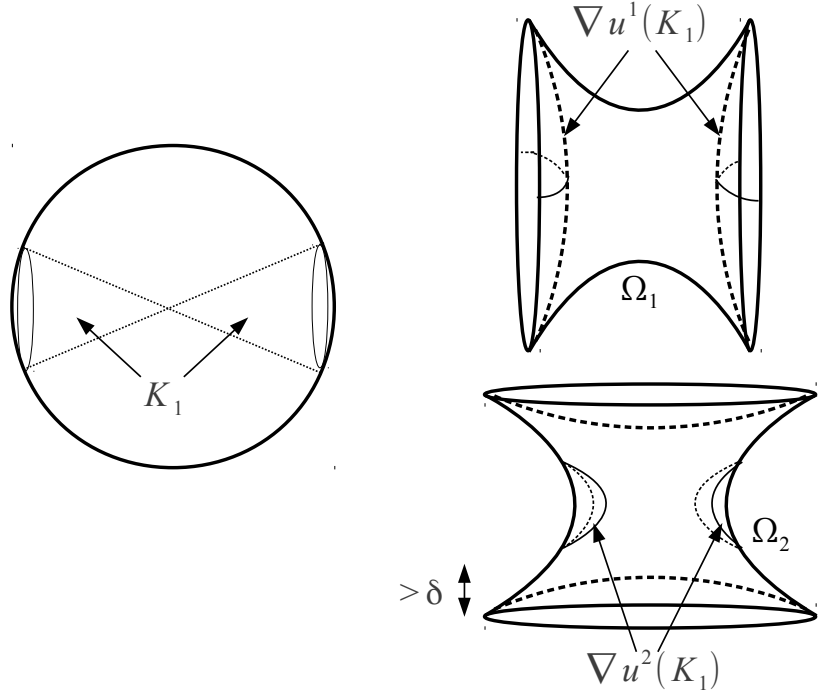}
 \caption{$\nabla u^1$ maps the cone $K_1$ to a region where $G_1$ is slightly non-convex, but $\nabla u^2$ maps it well
inside $\Omega_2$ where $G_2$ is uniformly convex.}
 \label{GradMap}
\end{figure}

Let $$\Omega_i = \nabla u^i(B_1 - K_i).$$ Note that $\Omega_1$ is just $\Omega_R$ squeezed by a factor of $1/2$ in the
$p_1$ direction. Let $\nu_i$ be the outer normals to $\Omega_i$. Since $u^i$ are homogeneous degree one we have
$\nu_i(\nabla u^i(x)) = x$ on $(B_1 - K_i) \cap S^2$. Furthermore, the preimage $x \in S^2$ of any point in $\Sigma$ satisfies
either $|x_1| \leq |x_3|$ or vice versa. It follows from these observations that if $(p^1, p^2) \in \Sigma$ then either
$$p^1 \in \Omega_1 \cap \{-\beta/2 \leq p^1_1 \leq \beta/2\} \text{   or   } p^2 \in \Omega_2 \cap \{-\beta/2 \leq p^2_3 \leq \beta/2\}$$
with $\beta$ such that $\varphi'(\beta)=1/2$, $\beta< 1-\delta$ (see figure \ref{GradMap}). Assume $p^1$ belongs to the set above.

Finally, let
$$F_0(p^1,p^2) = G_1(p^1) + G_2(p^2).$$
By rescaling Proposition \ref{KeyProposition} we have for $(p^1,p^2),\, (q^1,q^2) \in \Sigma$ that
\begin{align*}
S_{F_0}((p^1,p^2),(q^1,q^2)) &= S_{G_1}(p^1,q^1) + S_{G_2}(p^2,q^2) \\
& \geq \gamma|p^1-q^1|^2 - \epsilon |p^2-q^2|^2.
\end{align*}
Let $\omega_0 \in S^2$ be a preimage of $p^1$ under $\nabla u^1$. Then $|\nabla u^1(\omega)-\nabla u^1(\omega_0)| > c|\omega-\omega_0|$
and $|\nabla u^i(\omega)-\nabla u^i(\omega_0)| < C|\omega - \omega_0|$ for any $\omega \in S^2$, so
$$|p^2-q^2| \leq C|p^1-q^1|,$$
giving quadratic separation. By Corollary \ref{ExtensionCorollary} there is a smooth uniformly convex function $F$ on $\mathbb{R}^6$ so that $F = F_0$ and $ \nabla F = \nabla F_0$ on $\Sigma$, hence $\bf u$ satisfies the Euler-Lagrange system $\text{div}(\nabla F(D {\bf u})) = 0$ in $B_1 \setminus \{0\}$. Now it is straightforward to check that $\bf u$ is a weak solution of the system in the whole $B_1$. Indeed
$$\int_{B_1}\nabla F(D {\bf u})\cdot D  \psi=0, \quad \quad \forall  \psi \in C_0^\infty(B_1),$$
follows by integrating first by parts in $B_1 \setminus B_\epsilon$ and then letting $\epsilon \to 0$.
\section{Constructions}\label{Constructions}
In this section we prove the key step, Proposition \ref{KeyProposition}. Since $\Omega = \nabla u_0(B_1)$ is the surface obtained by revolving
$\Gamma$ around the $p_1$ axis, we can reduce to a one-dimensional problem on $\Gamma$ and then revolve the resulting picture
around the $p_1$ axis.
Since all of our constructions will be on $\mathbb{R}^2$ in this section we use coordinates $(x,y)$ rather than $(p_1,p_2)$.

\subsection{Setup}
Define $H$ to be an even function in $x$ and $y$ which has the form
\begin{equation}\label{Hdef}
H(x,y) = f(x) + h(x)(|y|-\varphi(x)),
\end{equation}
and is defined in a neighborhood of every point on $\Gamma_1 \cup \Gamma_3$, 
for some smooth functions $f$ and $h$ on $[-1,1]$. In our construction $h$ will be identically zero and $f$ linear near $x = \pm 1$,
so $H$ is linear in a neighborhood of the cusps of $\Gamma$. Notice that we can
extend $H$ to be a linear function (depending only on $x$) in a whole neighborhood of $\Gamma_2$ and similarly on $\Gamma_4$. Then $H$ is defined and smooth in a neighborhood of $\Gamma$.

\subsection{Inequalities for $\varphi$}
We now record some useful properties of $\Gamma$. For proofs see Section \ref{Appendix}. The first estimate gives an expansion for $\varphi$ near $x = -1$.

\begin{prop}\label{PhiComputations}
The function $\varphi$ is even, uniformly convex, and tangent to $y=|x|$ at $x = \pm 1$. Furthermore, $\varphi''$ is decreasing
near $x = -1$ and we have the expansion
\begin{equation}\label{phiexpansion}
\varphi''(-1 + \epsilon) = \sqrt{\frac{2}{3}} \epsilon^{-1/2} + O(1).
\end{equation}
\end{prop}

\begin{figure}
 \centering
    \includegraphics[scale=0.30]{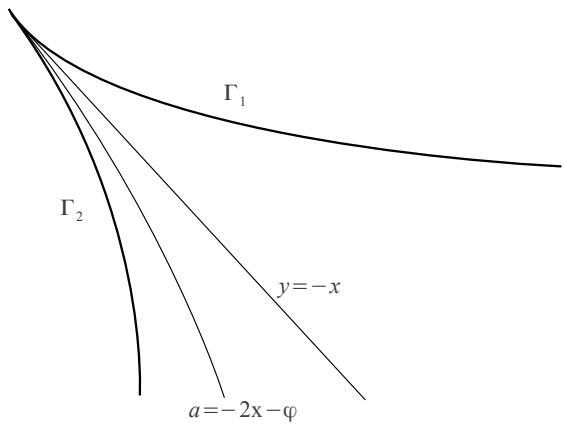}
 \caption{The graph $a = -2x - \varphi$ lies strictly above $\Gamma_2$.}
 \label{SideSepPic}
\end{figure}

The second estimate says that the vertical reflection of $\varphi$ over its tangent $y = -x$ lies above and separates from $\Gamma_2$ (see 
figure \ref{SideSepPic}). It follows easily from the uniform convexity of $\varphi$.

\begin{prop}\label{LeftSideSeparation}
The function $a(x) = -2x-\varphi$ is uniformly concave, tangent to $\Gamma_2$ at $x = -1$, and lies strictly above $\Gamma_2$ for $x > -1$. 
\end{prop}

\subsection{Euler-Lagrange Equation}
Let 
$$G(p_1,p_2,p_3) = H\left(p_1,\sqrt{p_2^2+p_3^2}\right).$$
The condition that $u_0$ solves the Euler-Lagrange equation $\text{div}(\nabla G(\nabla u_0)) = 0$ is equivalent to 
\begin{equation}\label{EulerLagrange}
h(x) = \frac{f''(x)}{2\varphi''(x)}.
\end{equation}

Indeed, since $G$ is linear near the surfaces obtained by revolving $\Gamma_2$ and $\Gamma_4$, we just need to verify the Euler-Lagrange equation
where $\nabla u_0$ is on the surface $\Omega_R$ obtained by revolving $\Gamma_1$. By passing a derivative the Euler-Lagrange equation $\text{div}(\nabla G(\nabla u_0))$
is equivalent to
$$\text{tr}\left(D^2G(\nabla u_0) \cdot D^2u_0\right) = 0.$$

Let $\Omega_R$ have outer normal $\nu$ and second fundamental form $II$.
Since $u_0$ is homogeneous degree one we have $\nu(\nabla u_0(x)) = x$ on $S^2$. Let $T$ be a frame tangent to $S^2$ at $x$, and differentiate to 
obtain $D_T^2u_0(x) = II^{-1}(\nabla u_0(x)).$
In coordinates tangent to $\Omega_R$ at $p = (p_1,\varphi(p_1),0)$ one computes
$$II = \frac{1}{\sqrt{1 + \varphi'^2}} \left( \begin{array}{cc}
\frac{\varphi''}{1 + \varphi'^2} & 0 \\
0 & -\frac{1}{\varphi} \end{array} \right), \quad
D^2G = \left( \begin{array}{cc}
\frac{f'' - h\varphi''}{1 + \varphi'^2} & 0 \\
0 &  \frac{h}{\varphi} \end{array} \right)$$
and the Euler-Lagrange formula follows.

\begin{rem}\label{ELComputation}
For a fast way to compute $D^2G$ in tangential coordinates, differentiate the equation $G(p_1,\varphi(p_1),0) = f(p_1)$:
$$\nabla G \cdot (1,\varphi') = f', \quad (1,\varphi')^T \cdot D^2G \cdot (1,\varphi') + h \varphi'' = f''.$$
The other eigenvalue comes from the rotational symmetry of $G$ around the $p_1$ axis.
\end{rem}

\begin{rem}\label{HigherDimEL}
If we do the computation in $\mathbb{R}^n$ we have $n-1$ rotational principal curvatures and derivatives, giving the Euler-Lagrange equation $h = \frac{f''}{(n-1)\varphi''}$.
\end{rem}

\subsection{Convexity Conditions}
Since most of our analysis is near a cusp, it is convenient to shift the picture by the vector $(1,-1)$ so that $\varphi, \,f$ are defined
on $[0,2]$ and $\varphi$ is tangent to $y = -x$ at zero. We assume this for the remainder of the section.

We examine convexity conditions between two points on $\Gamma_1$. Let $p = (x_0,\varphi(x_0))$ and $q = (x,\varphi(x))$. We first write
the equation for the tangent plane $L_p$ to $H$ at $p = (x_0,\varphi(x_0))$:
$$L_p(x,y) = f(x_0) + f'(x_0)(x-x_0) + h(x_0) \left[y- (\varphi(x_0) + \varphi'(x_0)(x-x_0)) \right].$$
Applying the Euler-Lagrange equation \eqref{EulerLagrange} we obtain
\begin{equation}\label{TangentPlane}
 L_p = f(x_0) + f'(x_0)(x-x_0) - \frac{f''(x_0)}{2\varphi''(x_0)} \left[y- (\varphi(x_0) + \varphi'(x_0)(x-x_0)) \right].
\end{equation}

By definition,
$$S_{H}(p,q) = f(x) - L_p(x,\varphi(x)).$$
Using equation \eqref{TangentPlane} we obtain
\begin{equation}\label{TopSeparation}
S_H(p,q) = \int_{x_0}^x f''(t)(x-t)\,dt - \frac{f''(x_0)}{2\varphi''(x_0)} \int_{x_0}^x \varphi''(t)(x-t)\,dt.
\end{equation}

\begin{defn}\label{Sep}
For a nonnegative function $g : \mathbb{R} \rightarrow \mathbb{R}$ define the weighted average
$$s_g(x_0,x) = \frac{\int_{x_0}^x g(t)(x-t)\,dt}{g(x_0)(x-x_0)^2}.$$
\end{defn}

With this definition we have
\begin{equation}\label{TS1}
S_H(p,q) = f''(x_0)  \left( s_{f''}(x_0,x) - \frac{1}{2}s_{\varphi''}(x_0,x) \right) (x-x_0)^2,
\end{equation}
thus, the first qualitative convexity condition is
\begin{equation}\label{FirstConvexityCondition}
s_{f''}(x_0,x) \geq \frac{1}{2}s_{\varphi''}(x_0,x).
\end{equation}

\begin{rem}\label{LocalComputationFirstCondition}
Notice that 
$$\lim_{x \to x_0} s_g(x_0,x)=\frac 12.$$
It is easy to check that if $g$ is increasing (decreasing) then $s_g(x_0,x)$ is increasing (decreasing) with $x$.
With this observation one verifies that condition \eqref{FirstConvexityCondition} holds for $x_0, \,x$ near $0$ if 
$f''(x) = Cx^{1-\alpha}$ for any $\alpha \in (0,1)$. Indeed, since $f''$ is increasing and $\varphi''$ is decreasing
one only needs to check the condition at $x = 0$, where one computes $s_{f''}(x_0,0) = \frac{1}{3-\alpha}$ and 
$\frac{1}{2}s_{\varphi''}(x_0,0) = \frac{1}{3} + O(\sqrt{x_0})$ which follows by Proposition \ref{PhiComputations}.
\end{rem}

We now examine convexity conditions between $p \in \Gamma_1$ and $q \in \Gamma_2$. 

Let $p = (x_0, \varphi(x_0))$.
In our construction we will have $h \geq 0$, and since $H$ is linear near $\Gamma_2$, we see that $S_H(p,q) \geq 0$ if the intersection line of tangent planes
to $H$ at $p$ and at $0$ lies above the line $y = -x$ on $[0,2]$. Using equation \eqref{TangentPlane} we compute the formula for the intersection line:
\begin{equation}\label{IntersectionLine}
\begin{split}
 y = \varphi(x_0) - \frac{2\varphi''(x_0)}{f''(x_0)}\int_{0}^{x_0} f''(t)(x_0-t)\,dt \\
 + \left(\varphi'(x_0) - \frac{2\varphi''(x_0)}{f''(x_0)}\int_{0}^{x_0}f''(t)\,dt\right)\cdot (x-x_0).
\end{split}
\end{equation}
If condition \eqref{FirstConvexityCondition} holds at $x = 0$, it means that the origin lies below the intersection line, thus $S_H(p,q) \geq 0$ for all $q \in \Gamma_2$ provided that the slope of the intersection line above is larger than $-1$:
$$ \varphi'(x_0) - \frac{2\varphi''(x_0)}{f''(x_0)}\int_{0}^{x_0}f''(t)\,dt \ge -1=\varphi'(0).$$

\begin{defn}\label{DerivDiff}
For a nonnegative function $g: \mathbb{R} \rightarrow \mathbb{R}$ define
$$d_g(x) = \frac{\int_0^{x} g(t)\,dt}{xg(x)}.$$
\end{defn}

With this definition the slope condition above can be written as
\begin{equation}\label{SecondConvexityCondition}
d_{f''}(x) \leq \frac{1}{2} d_{\varphi''}(x).
\end{equation}

\begin{rem}\label{LocalComputationSecondCondition}
Near $x = 0$ one computes $\frac{1}{2}d_{\varphi''}(x) = 1 + O(\sqrt{x})$. Thus, if $f''(x) = Cx^{1-\alpha}$ near $x = 0$ then \eqref{SecondConvexityCondition} holds.
However, away from a small neighborhood of $0$, condition \eqref{SecondConvexityCondition} will not hold in our construction. We will use formula \eqref{TangentPlane} 
more carefully, combined with Proposition \ref{LeftSideSeparation}, to deal with these cases.
\end{rem}

\begin{rem}\label{LinPartIndependence}
 Conditions \eqref{FirstConvexityCondition} and \eqref{SecondConvexityCondition} are independent of the linear part of $f$.
 Thus, when checking convexity conditions we only need to use the properties of $f''$.
\end{rem}

\subsection{Preliminary Construction}
As a stepping stone to proving Proposition \ref{KeyProposition} we construct first a $C^{1,\alpha}$ function $H_0$ near $\Gamma$, that is globally convex. 
We will use this construction to prove Theorem \ref{ScalarExample} in Section \ref{Appendix}. The function $H\in C^\infty$ is obtained by perturbing $H_0$. Below we define
$$G_0(p_1,p_2,p_3) = H_0\left(p_1,\sqrt{p_2^2 + p_3^3}\right).$$
Recall in the constructions below that we have shifted the picture by $(1,-1)$.

\begin{prop}\label{H0Construction}
 For any $\alpha \in (0,1)$ there exist a function $H_0$ near $\Gamma$ such that
 \begin{enumerate}
 \item $H_0$ is a linear function depending only on $x$ on $\Gamma_2$, and similarly on $\Gamma_4$.
  \item $H_0$ is pointwise $C^{1,1-\alpha}$ on the cusps of $\Gamma$ and smooth otherwise,
  \item $\text{div}(\nabla G_0(\nabla u_0)) = 0$ away from the cone $\{x_1^2 = x_2^2+x_3^2\},$
  \item $S_{H_0}(p,q) \geq 0$ for all $p,\,q \in \Gamma$,
  \item If $p = (x,\varphi(x))$ then $S_{H_0}(p,q) \ge \eta(x)|p-q|^2$ for all $q \in \Gamma$, where $\eta$ is 
  some continuous function on $[0,2]$ with $\eta > 0$ on $(0,2)$ and $\eta(0)=\eta(2)=0$.
 \end{enumerate}
\end{prop}

We will define $f_0$ by $f_0(0) = f_0'(0) = 0$ and prescribe $f_0''$, and then let $H_0$ be the function determined by $f_0$ 
through the Euler-Lagrange relation \eqref{EulerLagrange}. It is easy to check that condition \eqref{FirstConvexityCondition} holds if
we take $f_0'' = \varphi''$. However, we want $h_0 = f_0''/(2\varphi'')$ to go to zero at the endpoints so that
$H_0$ is linear on $\Gamma_2$ and $\Gamma_4$.

\begin{figure}
 \centering
    \includegraphics[scale=0.35]{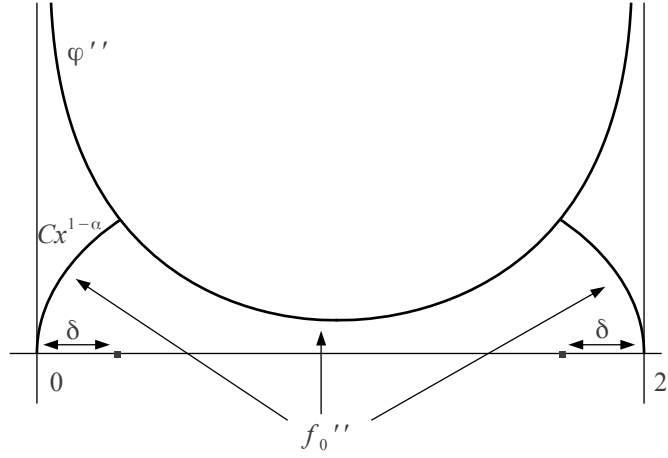}
 \caption{$f_0''$ agrees with $\varphi''$ on $[\delta,2-\delta]$, behaves like $x^{1-\alpha}$ near zero, and is symmetric
 around $x = 1$.}
 \label{f0pic}
\end{figure}

Motivated by the above and Remarks \ref{LocalComputationFirstCondition} and \ref{LocalComputationSecondCondition}, define
$$f_0''(x) = \begin{cases} \delta^{\alpha - 1}\varphi''(\delta)x^{1-\alpha}, \quad 0 \leq x \leq \delta \\
\varphi''(x), \quad \delta \leq x \leq 1 \\
f_0''(2-x), \quad 1 \leq x \leq 2
\end{cases}$$
(See figure \ref{f0pic}). Assume $\delta$ is tiny so that $\varphi''$ is well approximated by its expansion \eqref{phiexpansion}. 
Let $H_0$ be the function as in \eqref{Hdef} determined by $f_0$ through the Euler-Lagrange relation \eqref{EulerLagrange}.

\begin{proof}[{\bf Proof of proposition \ref{H0Construction}}]
The first three items are clear by construction so we check the convexity conditions. By symmetry we only need to consider $p \in \Gamma_1 \cup \Gamma_2$.
 
If $p \in \Gamma_2$ the positive separation is a consequence of $H_0 \geq 0$. This follows from the definition of $H_0$ on $\Gamma_1 \cup \Gamma_3$. 
Also, by symmetry, the linear function on $\Gamma_4$ intersects the linear function on $\Gamma_2$ on the vertical line $\{x=1 \}$, 
and since $\Gamma_4 \subset \{x>1\}$ we obtain $H_0 \ge 0$ on $\Gamma_4$ as well.

We now consider the situation when $p\in \Gamma_1$ and distinguish two cases depending 
whether $q \in \Gamma_1 \cup \Gamma_3$ or $q \in \Gamma_2 \cup \Gamma_4$. 

Let $p  = (x_0,\varphi(x_0))$.

{\bf First Case:} Assume first that $q = (x,\varphi(x)) \in \Gamma_1$. By symmetry of $f_0''$ around $x = 1$ we may assume $x < x_0$. 

If $x_0 \in [0,\delta]$ then by formula \eqref{TS1} and Remark \ref{LocalComputationFirstCondition} we have
$$S_{H_0}(p,q) \geq c(\alpha)f''(x_0)(x-x_0)^2.$$

If $x_0 \in [\delta, 2-\delta]$ we have $f_0'' = \varphi''$, so one computes
$$S_{H_0}(p,q) = \int_{x}^{x_0} (f''(t) - \frac 12 \varphi''(t))(t-x)\,dt.$$
If $x \geq \delta$ then this is clearly controlled below by $\frac{1}{4}\min(\varphi'')(x-x_0)^2$, and if $x < \delta$ then we have
$$S_{H_0}(p,q) = \int_{x}^{\delta} (f''(t) - \varphi''(t)/2)(t-x)\,dt + \frac{1}{2}\int_{\delta}^{x_0} \varphi''(t)(t-x)\,dt,$$
which is controlled below by 
$$\varphi''(\delta)(s_{f''}(\delta,x) - s_{\varphi''}(\delta,x)/2)(\delta - x)^2 + \frac{1}{4}\min(\varphi'')(x_0-\delta)^2 \geq c(\alpha)(x-x_0)^2.$$

Finally, if $x_0 \geq 2-\delta$ then since $f_0''/\varphi''$ is decreasing on $[\delta,2]$, we compute for $x \geq \delta$ that
$$S_{H_0}(p,q) \geq \frac{1}{2} \int_{x_0}^x f''(t)(x-t)\,dt \geq \frac{1}{4}\min\{f_0''(x_0),\min(\varphi'')\}(x-x_0)^2.$$
If $x < \delta$ then, since $f_0'' \leq \varphi''$ and they agree on $[\delta,2-\delta]$, we have using expansion \eqref{phiexpansion} that
$$S_{H_0}(p,q) \geq \frac{1}{2}\int_{\delta}^{2-\delta} \varphi''(t)(t-x)\,dt - C\sqrt{\delta} \geq c(x-x_0)^2.$$ 

If $q \in \Gamma_3$ then quadratic separation holds as well since 
$$\partial_yH_0(x_0,\varphi(x_0)) = \frac{f_0''(x_0)}{2\varphi''(x_0)} > 0.$$ 

{\bf Second Case:} By symmetry we may assume $q \in \Gamma_2$. If $x_0 \leq \delta$ we compute
$$d_{f_0''}(x_0) = \frac{1}{2-\alpha} < 1.$$
By Remark $\ref{LocalComputationSecondCondition}$, inequality \eqref{SecondConvexityCondition} holds strictly.

Now assume $x_0 \in [\delta,1]$. Define
$$g(x) = \varphi(x) - 2f_0(x).$$
Using the tangent plane formula \eqref{TangentPlane} we compute
$$L_p(x,g(x)) = -\int_{x_0}^x f''(t)(x-t)\,dt + \frac{1}{2} \int_{x_0}^x \varphi''(t)(x-t)\,dt = -S_{H_0}(p, (x,\varphi(x))) \leq 0$$
by the computations in the first case. Furthermore, since $f_0'' \leq \varphi''$, the graph of $g$ lies above the function 
$$a(x) = -2x-\varphi(x)$$
defined in Proposition \ref{LeftSideSeparation} (see figure \ref{SideSepKeyPic}). Since $a(x)$ lies strictly above $\Gamma_2$ for $x > 0$ and
$\partial_yH_0(x_0,\varphi(x_0)) = 1/2$, we have strictly positive separation on $\Gamma_2$.

\begin{figure}
 \centering
    \includegraphics[scale=0.35]{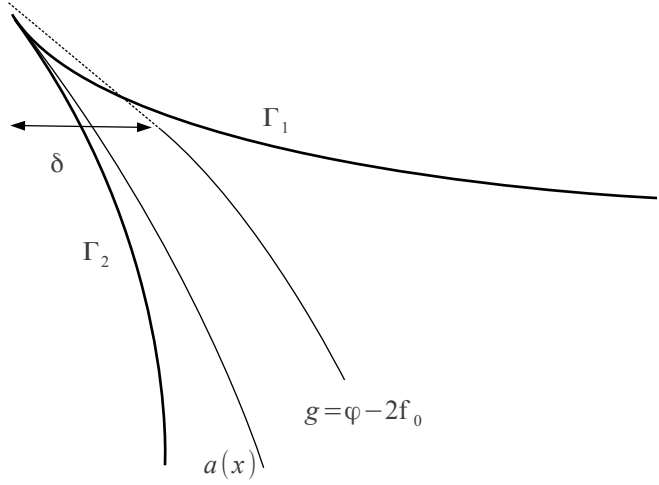}
 \caption{The tangent plane at $(x,\varphi(x))$ is negative on the curve $y = g(x)$, hence on $\Gamma_2$, for $x \in [\delta,2-\delta]$.}
 \label{SideSepKeyPic}
\end{figure}

Finally, for $x_0 \in [1, 2]$, the intersection of the tangent planes at $p$ and at $\tilde{p} = (2-x_0,\varphi(2-x_0))$ is the line $x = 1$ since
$f_0''$ is symmetric around $x = 1$. By the previous computations, the tangent plane at $\tilde{p}$ is negative on $\Gamma_2$. Thus, the tangent
plane at $p$ is negative on $\Gamma_2$, completing the proof.
\end{proof}

\subsection{Proof of Key Proposition}
We can slightly modify the construction of $H_0$ from the previous section to make it smooth, at the expense of giving up a little convexity near the cusps of $\Gamma$.
Below $\delta, \, \gamma > 0$ are small constants depending only on $\varphi$.
Let $G(p_1,p_2,p_3) = H\left(p_1,\sqrt{p_2^2+p_3^2}\right)$.

\begin{prop}\label{HConstruction}
 For any $\epsilon > 0$ there exists a smooth function $H$ defined on a neighborhood of $\Gamma$ such that
 \begin{enumerate}
 \item $H$ is linear (depending only on $x$) in a neighborhood of $\Gamma_2$, respectively $\Gamma_4$,
  \item $\text{div}(\nabla G(\nabla u_0)) = 0$,
  \item $H_y(x,\varphi(x)) \geq \frac 1 2$ for $x \in [\delta, 2-\delta]$, and $H_y \geq 0$ on $\Gamma_1$,
  \item If $p = (x,\varphi(x))$ with $x \in [\delta,2-\delta]$ then $S_{H}(p,q) \geq \gamma |p-q|^2$ for all $q \in \Gamma$,
  \item $S_{H}(p,q) \geq - \epsilon |p-q|^2$ otherwise for $p,\,q \in \Gamma$.
 \end{enumerate}
\end{prop}

Note that the key Proposition \ref{KeyProposition} follows easily from Proposition \ref{HConstruction} by defining $G$ as above.

\begin{figure}
 \centering
    \includegraphics[scale=0.35]{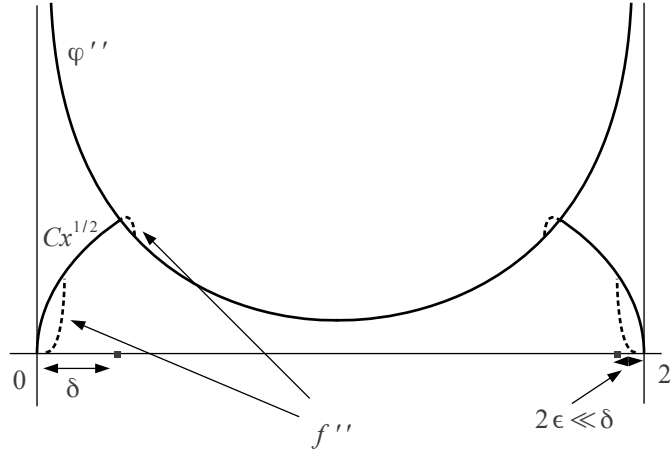}
 \caption{$f''$ is a small perturbation of $f_0''$ that connects smoothly to $\varphi''$ near $x = \delta$ and goes quickly to zero near $x = 0$.}
 \label{fpic}
\end{figure}

Let $\alpha = \frac{1}{2}$ in the construction of $f_0''$ from the previous section and let $\epsilon \ll \delta$. Let $f''$ be a smoothing of $f_0''$ defined by cutting it off smoothly to zero between
$\epsilon$ and $2\epsilon$, gluing it smoothly to itself between $\delta$ and $\delta + \epsilon$, and making it symmetric over $x = 1$ (see figure \ref{fpic}). 
Let $H$ be the function in \eqref{Hdef} determined by $f$ through
the Euler-Lagrange relation \eqref{EulerLagrange}.

\begin{figure}
 \centering
    \includegraphics[scale=0.35]{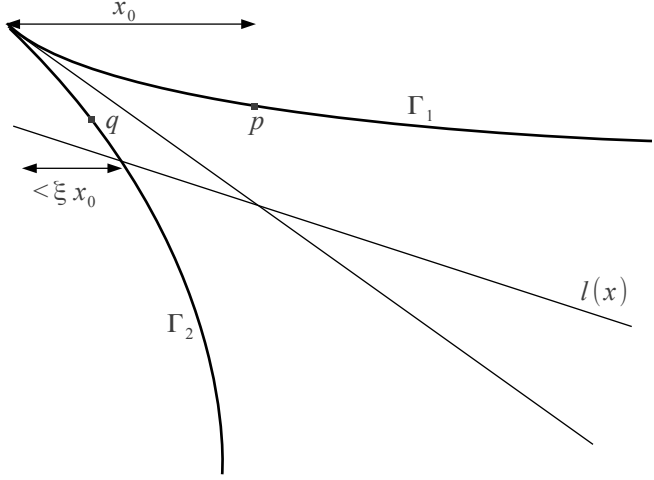}
 \caption{The separation is positive if $q$ is below the line $l(x)$, and if the separation is negative then $|p-q|$ is of order $x_0$.}
 \label{DistSepPic}
\end{figure}

\begin{proof}[{\bf Proof of proposition \ref{HConstruction}}]
The first three conclusions are clear by construction so we just need to check the convexity conclusions. Most of them will follow by continuity.

If $p \in \Gamma_2$ we have positive separation since $H \geq 0$, so assume $p = (x_0,\varphi(x_0))$.

If $x_0 \in [\delta,2-\delta]$ then the conclusion holds by continuity from the arguments
in the proof of Proposition \ref{H0Construction} after taking $\epsilon$ small. 

Next we may assume by symmetry that $x_0 \in [0, \delta]$.

{\bf Case 1:} Assume that $x_0 \geq 10\epsilon$. 

If $q = (x,\varphi(x))$ with $x > x_0$ then the positive separation follows again by continuity.
If $x < x_0$ one computes
$$2s_{f''}(x_0,x) \geq 2s_{f''}(x_0,0) \geq \frac{4}{5}(1-(1/5)^{5/2}) > s_{\varphi''}(x_0,0) \geq s_{\varphi''}(x_0,x)$$
so condition \eqref{FirstConvexityCondition} holds and we have positive separation on $\Gamma_1$. 

Since the cutoff is between $\epsilon$ and $2\epsilon$ and $f''$ is increasing for $x < \delta$ we compute
\begin{equation}\label{LeftSepKey}
d_{f''}(x_0) < \frac{2}{3}.
\end{equation}
and by Remark \ref{LocalComputationSecondCondition} the condition \eqref{SecondConvexityCondition} holds for $x_0 < \delta$.
We thus have positive separation on $\Gamma_2$ and $\Gamma_3$.

Finally, for $q \in \Gamma_4$ positive separation follows again by continuity.

This establishes positive separation everywhere for $x_0 \in [10\epsilon, 2-10\epsilon]$.

{\bf Case 2:} Assume $x_0 \leq 10\epsilon$. 

The tangent plane at $p$ is of order $\epsilon$ on $\Gamma$, so we have
positive separation when $q \in \Gamma_4$. 

Using that $f''$ is increasing and $\varphi''$ decreasing near $0$, we obtain positive separation if $q = (x,\varphi(x))$ with $x \in [x_0, \delta]$. The same holds for $x> \delta$ by continuity.

If $q = (x,\varphi(x))$ for $x < x_0$ we compute
$$S_{H}(p,q) \geq - f''(x_0)s_{\varphi''}(x_0,x)(x-x_0)^2 \geq - C\sqrt{\epsilon} \, \,|p-q|^2,$$
since $s_{\varphi''}(x_0,x) \le s_{\varphi''}(x_0,0) \le 1$. This gives the desired estimate on $\Gamma_1$. 

Next we bound $S_H(p,q)$ with $q \in \Gamma_2$. For this we estimate the location of the intersection line $l_p$ of the tangent plane at $p$ with $0$.
By \eqref{IntersectionLine}, $l_p$ passes through
$$\left(x_0, \varphi(x_0) - \frac{2\varphi''(x_0)}{f''(x_0)}\int_{0}^{x_0} (x_0-t)f''(t)\,dt\right).$$
We first claim that this point lies above the line $y = -x$. Indeed,
since $f''$ is increasing in $[0,x_0]$, the second component is larger than $\varphi(x_0) - \varphi''(x_0)x_0^2$, and using the expansion \eqref{phiexpansion} we see that
$$\varphi(x_0) + x_0 \geq \left(\frac{4}{3}\varphi''(x_0) + O(1)\right)x_0^2 > \varphi''(x_0)x_0^2.$$

By \eqref{LeftSepKey} the slope of $l_p$ is between $-1$ and $0$, so we have positive separation for $q \in \Gamma_3$ and
$q \in \Gamma_2 \cap \{y < -x_0\}$.

Finally, from \eqref{IntersectionLine} we see that the slope of $l_p$ is less than $\varphi'(x_0)$. Thus,
for $x < x_0$, $l_p$ lies above the line
$$y=l(x) = -x_0 + \varphi'(x_0)(x-x_0).$$
A short computation using the expansion \eqref{phiexpansion} shows that $l(x)$ crosses $a(x)$, hence $\Gamma_2$, at some $x < \xi x_0$ where
$$\xi + \frac{2}{3}\xi^{3/2} = 1 + O(\sqrt{\epsilon}).$$
In particular, $\xi < 1-c$. This gives that the separation is positive on $\Gamma_2 \cap \{x > \xi x_0\}$, and otherwise the separation is 
at worst $-C\sqrt{\epsilon} x_0^2 \geq -C\sqrt{\epsilon}|p-q|^2$ (see figure \ref{DistSepPic}).

\end{proof}

\begin{rem}\label{ConvexityLoss}
 The proof shows in fact that $S_{H}(p,q)$ is only negative for $p,\,q$ very close to the same cusp.
\end{rem}

\section{Appendix}\label{Appendix}
\subsection{Convex Extension Lemma}

\begin{proof}[{\bf Proof of lemma \ref{ExtensionLemma}}]
Let ${\bf v}_{T}$ be the tangential component and ${\bf v}_{\perp}$ be the normal component, and let $\nabla^{\Sigma}G$ be the gradient of $G$ on $\Sigma$. Note that condition \ref{QuadSepCondition} 
implies ${\bf v}_T = \nabla^{\Sigma}G$. For $x \in \Sigma$ let $T_x, \, N_x$ be the tangent  and normal subspaces to $\Sigma$ at $x$. Let $d_{\Sigma}(y)$ be the distance from $y$ to $\Sigma$ and let 
$$\Sigma^r = \{y: d_{\Sigma}(y) < r\}.$$
Finally, for $x \in \mathbb{R}^n$ let $y(x)$ be the closest point in $\Sigma$ to $x$. 
It is well-known that $y$ and $d_{\Sigma}^2$ are well-defined and smooth in a neighborhood of $\Sigma$, and for $x \in \Sigma$, $D_xy(x)$ is the projection to $T_x$ and $D^2(d_{\Sigma}^2/2)(x)$ is the 
projection to $N_x$. (For proofs, see for example \cite{AS}).

{\bf Step 1:} We claim that the function
$$F(x) = G(y(x)) + {\bf v}(y(x)) \cdot (x-y(x)) + \frac{A}{2}d_{\Sigma}^2(x)$$
with $A$ large lifts quadratically from its tangent planes in $\Sigma^{\sigma}$ for $\sigma$ sufficiently small. We first compute for $x \in \Sigma$ that
\begin{align*}
F(x + \epsilon z) &= G(x) + \epsilon\left(\nabla^{\Sigma}G(x) \cdot z_{T} + {\bf v}(x) \cdot z_{\perp}\right) + O(\epsilon^2) \\
&= G(x) + \epsilon {\bf v}(x) \cdot z + O(\epsilon^2)
\end{align*}
giving that $F = G$ on $\Sigma$ and $\nabla F = {\bf v}$ on $\Sigma$.

Now, for $\epsilon$ small and $\nu \in N_x$ we have $y(x + \epsilon \nu) = x$ and $d_{\Sigma}(x+\epsilon \nu) = \epsilon$, so $F_{\nu\nu}(x) = A.$
In addition,  if $x \in \Sigma$ and $x + \epsilon z \in \Sigma$ for some unit vector $z$ then by hypothesis we have
\begin{align*}
F(x + \epsilon z) &= F(x) + \epsilon \nabla F(x) \cdot z + \frac{\epsilon^2}{2}z^T \cdot D^2F(x) \cdot z + O(\epsilon^3) \\ 
&\geq F(x) + \epsilon \nabla F(x) \cdot z + \gamma \epsilon^2.
\end{align*}
Taking $\epsilon$ to zero we see that $F_{\tau\tau}(x) > 2\gamma$ for any tangential unit vector $\tau$.

Take any unit vector $e$ and write $e = \alpha \tau + \sqrt{1-\alpha^2} \nu$ for some unit $\tau \in T_x\Sigma$ and $\nu \in N_x\Sigma$. Since $D^2(d_{\Sigma}^2/2)$ is the projection
matrix onto $N_x$ at $x \in \Sigma$, we have
\begin{align*}
F_{ee}(x) &= \alpha^2F_{\tau\tau} + (1-\alpha^2)F_{\nu\nu} + 2\alpha\sqrt{1-\alpha^2}(F - Ad_{\Sigma}^2/2)_{\tau\nu} \\
&\geq 2\alpha^2\gamma + (1-\alpha^2)A - C\alpha\sqrt{1-\alpha^2}
\end{align*}
for some $C$ independent of $A$. We conclude that $D^2F > \frac{3}{2}\gamma I$
on $\Sigma$ for $A$ sufficiently large, and in particular, $D^2F > \frac{3}{2}\gamma I$ on a neighborhood $\Sigma^{2\rho}$ of $\Sigma$.

Finally, we show that the tangent planes to $F$ in $\Sigma^{\sigma}$ separate quadratically for $\sigma$ small. 
Let $x, \,z \in \Sigma^{\sigma}$. We divide into two cases.

If $|z-x| < \rho$ then $x$ and $z$ can be connected by a line segment contained in $\Sigma^{2\rho}$, so it is clear that 
$$F(z) > F(x) + \nabla F(x) \cdot (z-x) + \frac{3}{4}\gamma|z-x|^2.$$

If on the other hand $|z-x| > \rho$, we use that
$$F(y(z)) > F(y(x)) + \nabla F(y(x)) \cdot (y(z)-y(x)) + \gamma |y(z)-y(x)|^2.$$
Replacing $y(z)$ by $z$ and $y(x)$ by $x$ changes these quantities by at most $C\sigma$, and since
and $|z-x| > \rho$ we have that
$$F(z) > F(x) + \nabla F(x) \cdot (z-x) + \frac{3}{4}\gamma |z-x|^2$$
for all $x,\,z \in \Sigma^{\sigma}$ for $\sigma$ small.

{\bf Step 2:} From now on denote the open set $\Sigma^{\sigma}$ by $N$.
Let $N_{\epsilon}$ denote $\{x \in N: B_{\epsilon}(x) \subset N\}$.
Finally, let $\rho_{\epsilon}$ denote the standard mollifier $\epsilon^{-n}\rho(x/\epsilon)$ where $\rho$ is supported in $B_1$, nonnegative, smooth 
and has unit mass.

We define a global uniformly convex function that agrees with $F$ on $N$. Let
$$H_0(y) = \sup_{x \in N}\left\{F(x) + \nabla F(x) \cdot (y-x) + \frac{3}{4}\gamma |y-x|^2\right\}.$$
Then $H_0$ is a uniformly convex function on $\mathbb{R}^n$ with $D^2H_0 \geq \frac{3}{2}\gamma I$ and furthermore  by construction we have that $H_0 = F$ on $N$.

To finish we glue $H_0$ to a mollification. Fix $\delta$ so that $\Sigma \subset N_{2\delta}$. Let 
$$H_{\epsilon} = \rho_{\epsilon} \ast H_0$$
for some $\epsilon$ small. In $N_{\delta}$ we have
$$|H_{\epsilon} - H_0|, \, |\nabla H_{\epsilon}-\nabla H_0| < C\epsilon.$$
Finally, since $D^2H_0 \geq \frac{3}{2}\gamma I$ we have $D^2H_{\epsilon} > \frac{3}{2}\gamma I.$

Let $\eta$ be a smooth cutoff function which is $1$ on $N_{2\delta}$ and $0$ outside of $N_{\delta}$. Then let
$$H= \eta H_0 + (1-\eta) H_{\epsilon}.$$
We compute
$$D^2H = \eta D^2H_0 + (1-\eta) D^2H_{\epsilon} + 2\nabla \eta \otimes \nabla(H_0-H_{\epsilon}) + D^2\eta(H_0-H_{\epsilon}).$$
Then $H$ is smooth, $H = F$ on $N_{2\delta}$ and taking $\epsilon$ small we have $D^2H > \gamma I$, completing the construction.
\end{proof}

\subsection{Expansion of $\varphi$}
\begin{proof}[{\bf Proof of proposition \ref{PhiComputations}}]
The symmetries of $\varphi$ follow from the symmetries of $w$. 

The curve $\Gamma_1$ is parametrized by $\nabla w(\theta)$ for $\theta \in [\pi/4, \,3\pi/4]$. Let $\nu$ be the upward normal to $\Gamma_1$.
Since $w$ is homogeneous degree one we have $\nu(\nabla w(\theta)) = \theta$. Differentiating we
get the the curvature $\kappa = \frac{1}{g'' + g}$ where $g(\theta) = \frac{-1}{\sqrt{2}}\cos 2\theta$ are the values of $w$ on $S^1$. 
Thus, $\varphi$ is uniformly convex and its second derivatives blow up near $x = \pm 1$. To quantify this we compute
\begin{align*}
 \nabla w(\theta) &= g(\theta)(\cos \theta,\, \sin \theta) + g'(\theta)(-\sin\theta,\, \cos \theta) \\
 &= \frac{1}{\sqrt{2}}(-\cos \theta(1+2\sin^2\theta),\, \sin \theta(1+2\cos^2\theta)).
\end{align*}
Expanding around $\theta = \frac{\pi}{4}$ (which gets mapped to the left cusp on $\Gamma_1$) we get
\begin{equation}\label{GradientExpansion}
\varphi\left(-1 + \frac{3}{2}\theta^2 + \theta^3 + O(\theta^4)\right) = 1-\frac{3}{2}\theta^2 + \theta^3 + O(\theta^4).
\end{equation}
Differentiating implicitly one computes
$$\varphi''(-1 + \epsilon) = \sqrt{\frac{2}{3}}\epsilon^{-1/2} + O(1)$$
and that $\varphi''$ is decreasing near $-1$.
\end{proof}

\subsection{Theorem \ref{ScalarExample}}

In \cite{DS} the authors show that if $u$ is a scalar minimizer to a convex functional $\int_{B_1} F(\nabla u)\,dx$ on $\mathbb{R}^2$ and $F$ is uniformly convex in a neighborhood of $\nabla u(B_1) \cap \{|p_1| < 1\}$
then $\nabla u$ cannot jump arbitrarily fast across the strip. In particular, $\nabla u(B_{\gamma})$ localizes to $\{p_1 < 1\}$ or $\{p_1 > -1\}$ for some $\gamma$
small. In this final section we use the preliminary construction $H_0$ from section \ref{Constructions} to indicate
why this result is not true in three or higher dimensions.

Make a global extension of $H_0$ by taking 
$$\bar H_0(x) = \sup_{p \in \Gamma_1 \cup \Gamma_3}\{H_0(p) + \nabla H_0(p) \cdot (x-p) + \eta(p_1)|x-p|^2\}.$$
The resulting extension is smooth near any non-cusp point of $\Gamma$. It is uniformly convex near each point on $(\Gamma_1 \cup \Gamma_3) \cap \{|p_1| < 1\}$ with the modulus of convexity decaying towards the cusps.
Furthermore, $\bar H_0$ is flat in a neighborhood of every point on $(\Gamma_2 \cup \Gamma_4) \cap \{|p_2| < 1\}$. Finally, if $p$ is a cusp of $\Gamma$ 
then it is straightforward to check that $\bar H_0$ is pointwise $C^{1,1-\alpha}$ at $p$, i.e. $S_{\bar H_0}(p, x) < C|x-p|^{2-\alpha}$ for all $x$ near $p$.
By iterating a mollification and gluing procedure similar to those used in the proof of lemma \ref{ExtensionLemma} near the cusps we can get a global convex 
extension $\bar H$ that is smooth away from the cusps, uniformly convex on $\Gamma_1 \cup \Gamma_3$ away from the cusps, flat on convex sets
containing $\Gamma_2$ and $\Gamma_4$, and $C^{1,1-\alpha}$ at the cusps.

\begin{rem}\label{DimensionReg}
In dimension $n$ the Euler-Lagrange equation allows us to take $f_0''(x) = x^{n-2-\alpha}$ near the cusp, which gives $\bar H$ an extra
derivative for each dimension.
\end{rem}

\begin{figure}
 \centering
    \includegraphics[scale=0.40]{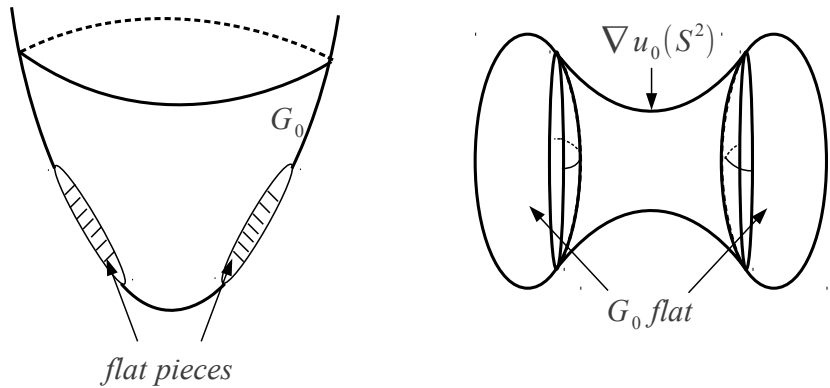}
 \caption{$G_0$ is linear on two bounded convex sets containing $\nabla u_0(\{|x_1| > r\})$.}
 \label{G0Pics}
\end{figure}

Let $G_0$ be the function on $\mathbb{R}^3$ obtained by revolving $\bar H$ around the $p_1$ axis (see figure \ref{G0Pics}).
By construction $u_0$ solves the Euler-Lagrange equation $\text{div}(\nabla G_0(\nabla u_0)) = 0$ away from the cone $C_0 = \{|x_1| = r\}$ where $r = \sqrt{x_2^2 + x_3^3}$.
Thus, it is not immediate that $u_0$ minimizes $\int_{B_1} G_0(\nabla u_0)\,dx$. However, we claim $u_0$ is a minimizer. To show this we must establish
$$\int_{B_1} \nabla G_0(\nabla u_0) \cdot \nabla \psi \,dx = 0$$
for any $\psi \in C^{\infty}_0(B_1)$. The contribution from integrating in $B_{\epsilon}$ and a thin cone $\{(1-\epsilon)r < |x_1| < (1+\epsilon)r\}$
is small. Integrating by parts in the remaining region with boundary $S$, we get a boundary term of the form $\int_S \psi \nabla G_0(\nabla u_0) \cdot \nu \,ds$
where $\nu$ is the outer normal. The cones $\{|x_1| = (1\pm \epsilon)r\}$ are $\epsilon$ close, and the outward normals on these cones are $\epsilon$ close
to flipping direction, so by the continuity of $\nabla G_0$ the contribution from this term is also small. Taking $\epsilon$ to zero we get the desired result.



 \section*{Acknowledgement}
 C. Mooney was supported by NSF fellowship DGE 1144155. 

O. Savin was supported by NSF grant DMS-1200701.


\frenchspacing
\bibliographystyle{plain}

\end{document}